\documentclass[12pt]{article}
\usepackage{amsmath}
\newtheorem{theorem}{Theorem}
\newcommand{\eqdef}{\, =\kern -12.7pt\raise 6pt\hbox{{\tiny\textrm{def}}}\,\,}
\def\hypf#1#2#3#4#5{\, _{#1}\kern -1pt F_{#2}\kern -3pt\left[
\begin{matrix}#3\\#4\end{matrix}\,;\,#5\right]}
\def\bars{\atopwithdelims||}
\def\gamp{\gamma\kern1pt '}

\begin{document}
\title{The method of characteristics, and `problem 89' of Graham, Knuth and Patashnik}
\author{Herbert S. Wilf\\
Department of Mathematics, University of Pennsylvania\\
Philadelphia, PA 19104-6395\\
\texttt{<wilf@math.upenn.edu>}
\\\\\\
To Richard Stanley, on his sixtieth birthday}
\date{}
\maketitle

\vspace{.8in}

\begin{abstract}
We apply the method of characteristics for the solution of pde's to two combinatorial problems. The first is finding an explicit form for a distribution that arises in bio-informatics. The second is a question raised by Graham, Knuth and Patashnik about a sequence of generalized binomial coefficients. We find an exact formula, which factors in an interesting way, in the case where one of the six parameters of the problem vanishes. We also show that the associated polynomial sequence has real zeros only, provided that one parameter vanishes, and the other five are nonnegative.
\end{abstract}

\vspace{1in}

\newpage

The method of characteristics is well known in the theory of partial differential equations, and it is widely used for solving linear pde's of the first order. An introductory description of the method and a number of examples are in \cite{fbh}. In combinatorics the method has scarcely been used, though perhaps it deserves wider recognition. The main idea is that a linear recurrence in more than one variable will lead to a linear partial differential equation for the generating function of the solution of that recurrence, and the method of characteristics might well be able to yield the solution. We give two examples below.

The first example is intended to introduce the method in a combinatorial context, and it concerns a rather special problem that was motivated by DNA sequence matching. The second example is drawn from a Research Problem in Graham, Knuth and Patashnik \cite{gkp}, concerning a two-variable recurrence that generalizes the familiar Pascal triangle and thereby defines a family of generalized binomial coefficients. In the first example we obtain a complete answer to the question, using the method of characteristics. In the second example we obtain a complete solution also, though its form is very unwieldy. In an important special case, though, we can obtain the solution in a nice form.

\section{A problem from bio-informatics}
The following question has application to the problem of deciding if two DNA samples match sufficiently. It was asked by Warren Ewens, of the University of Pennsylvania.

Let $Q_0(x)=1-x^n$ and
\begin{equation}
\label{eq:recurr}
Q_{k+1}(x)=Q_k(x)+\frac{1-x}{k+1}Q_k'(x)\qquad (k=0,1,2,\dots ).
\end{equation}
It is required to find a nice form for $Q_k$ and to study its properties.

We will in fact show three nice forms for $Q_k(x)$. The first one is
\begin{equation}
\label{eq:form1}
Q_k(x)=1-\sum_{j=0}^k{n\choose j}(1-x)^jx^{n-j},\qquad (k=0,1,2,\dots)
\end{equation}
Next we have, for $k=0,1,2,\dots$,
\begin{equation}
\label{eq:form2}
Q_k(x)=-(n-k){n\choose k}(-1)^{n+k}\sum_{r\ge 1}{k\choose n-r}\frac{(-1)^r}{r}x^r+
\begin{cases}
1,&\mathrm{if}\ k<n;\\
0,& \mathrm{else.}\\
\end{cases}
\end{equation}
But the nicest form is
\begin{equation}
\label{eq:form3}
Q_k(x)=(1-x)^{k+1}\sum_{j=0}^{n-k-1}{j+k\choose j}x^j.
\end{equation}
From this last form it is clear that
\[\lim_{n\to\infty}Q_k(x)=1\qquad (|x|<1)\]
for every $k\ge 0$. Thus for fixed $k\ge 0,|x|<1$, $Q_k(x)$ is a cumulative probability distribution function of $n=0,1,2,3,\dots $. Indeed, keeping $k$ and $x$ fixed, the corresponding probability density function is
\[p_n={n\choose k}x^{n-k}(1-x)^{k+1}\qquad (n=k, k+1, k+2,\dots ).\]

To prove these things, put $Q(x,t)=\sum_{k\ge 0}Q_k(x)t^k$. Then the recurrence (\ref{eq:recurr}) becomes, after the usual generating function operations,
\begin{equation}
\label{eq:pde1}
(x-1)\frac{\partial Q}{\partial x}+(1-t)\frac{\partial Q}{\partial t}=Q.
\end{equation}
This linear partial differential equation of the first order can be solved by the method of characteristics. To do that we find two independent solutions of the ordinary differential equations
\[\frac{dx}{x-1}=\frac{dt}{1-t}=\frac{dQ}{Q}.\]
We easily find such solutions in the form
\[(1-t)(1-x)=c_1;\qquad (1-t)Q=c_2.\]
Hence the general solution of the p.d.e. (\ref{eq:pde1}) is $c_2=f(c_1)$, where $f$ is an arbitrary function, i.e.,
$(1-t)Q(x,t)=f((1-t)(1-x))$. The initial data $Q(x,0)=1-x^n$ yields $f(u)=1-(1-u)^n$, hence
\[Q(x,t)=\frac{1-(1-(1-t)(1-x))^n}{1-t}.\]
After expanding the power above, the form (\ref{eq:form1}) follows. To get (\ref{eq:form2}) we expand the power of $1-x$ in (\ref{eq:form1}), collect powers of $x$, and observe that the coefficient of each power of $x$ is a Vandermonde convolution sum that can be evaluated in closed form.

Finally, we can verify form (\ref{eq:form3}) as follows. Expand the factor $(1-x)^{k+1}$ in (\ref{eq:form3}), to get
\begin{eqnarray*}
Q_k(x)&=&\sum_{\ell}(-1)^{\ell}{k+1\choose \ell}x^{\ell}\sum_{j\le n-k-1}{j+k\choose j}x^j\\
&=&\sum_r(-1)^rx^r\sum_{j=0}^{n-k-1}(-1)^j{k+1\choose r-j}{k+j\choose j}.
\end{eqnarray*}
But by Gosper's algorithm the last summand is
\[(-1)^j{k+1\choose r-j}{k+j\choose j}=g_{j+1}-g_j,\]
where
\[g_j=\frac{(-1)^{j+1}(k+j)!}{r(j-1)!k!}{k\choose r-j}.\]
Hence the inner sum is
\[g_{n-k}=\frac{(-1)^{n-k+1}n!}{r(n-k-1)!k!}{k\choose r-n+k}=(-1)^{n-k+1}\frac{n-k}{r}{n\choose k}{k\choose r-n+k},\]
and we have arrived back at form (\ref{eq:form2}) again, which completes the proof.

A glance at form (\ref{eq:form1}) also shows that $Q_k(x)$ vanishes for all $k\ge n$. For example, when $n=5$ the sequence $\{Q_k(x)\}_0^{\infty}$ is
\begin{eqnarray*}
&&1 - {x^5},1 - 5\,{x^4} + 4\,{x^5},1 - 10\,{x^3} + 15\,{x^4} - 6\,{x^5},
  1 - 10\,{x^2} + 20\,{x^3} - 15\,{x^4} + 4\,{x^5},\\
&&\qquad\qquad 1 - 5\,x + 10\,{x^2} - 10\,{x^3} + 5\,{x^4} - {x^5},0,0,0,\dots
\end{eqnarray*}

\section{The problem of Graham et al}
In \cite{gkp}, ``Research problem'' 89 of Chapter 6 asks for a general theory of recurrences of the form
\begin{equation}
\label{eq:rec}
{n\bars k}=(\alpha n+\beta k+\gamma){n-1\bars k}+(\alpha' n+\beta' k+\gamp){n-1\bars k-1}+\delta_{n,0}\delta_{k,0},\qquad (n,k\ge 0)
\end{equation}
which, of course, includes the recurrences for the binomial coefficients, the Stirling numbers of both kinds, the Eulerian numbers, and a variety of other     combinatorial sequences.

A number of contributions to this problem have been made by P. Th\'eor\^ et, in his dissertation \cite{th1} and in two subsequent papers \cite{th2}, \cite{th3}. In \cite{th1} he considered both eq. (\ref{eq:rec}), whose solutions he termed ``special hyperbinomial sequences,'' as well as the generalization in which the coefficients of (\ref{eq:rec}) are not required to be linear functions of $n$ and $k$, which were called ``hyperbinomial sequences.''

For technical reasons, we'll begin by rewriting their recurrence just slightly, as
\begin{equation}
\label{eq:recur}
{n+1\bars k}=(\alpha n+\beta k+\gamma){n\bars k}+(\alpha' n+\beta' k+\gamp ){n\bars k-1}.\qquad (n,k\ge 0; {0\bars k}=\delta_{k,0})
\end{equation}
\subsection{The generating function}
Following the usual generatingfunctionological techniques, we introduce the gf
\[u(x,y)=\sum_{n,k\ge 0}{n\bars k}\frac{x^n}{n!}y^k.\]
It is then easy to translate the recurrence (\ref{eq:recur}) into a linear partial differential equation of first order, viz.
\begin{equation}
\label{eq:pde}
(1-\alpha x-\alpha' y)\frac{\partial u(x,y)}{\partial x}=y(\beta +\beta' y)\frac{\partial u(x,y)}{\partial y}+(\gamma+(\beta' +\gamp )y)u(x,y),
\end{equation}
in $u(x,y)$, together with the initial condition $u(0,y)=1$. Such equations can in principle always be solved by the method of characteristics, though that method may require some tricks to make it work. Fortunately the Maple package PDEtools can save us the trouble of dealing with these tricks, and we can easily check that the results that it obtains are correct. We'll now do a few specific examples and then exhibit the general solution.
\begin{itemize}
\item \textit{Binomial coefficients}. With $(\alpha,\beta,\gamma,\alpha',\beta',\gamp )=(0,0,1,0,0,1)$, the \texttt{pdsolve} command produces the general solution
\[u(x,y)=F(y)e^{x(1+y)},\]
where $F$ is an arbitrary function. Since $u(0,y)=1=F(y)$, we have $u(x,y)=e^{x(1+y)}$.
\item \textit{Stirling numbers of the first kind}. Now $(\alpha,\beta,\gamma,\alpha',\beta',\gamp )=(1,0,0,0,0,1)$ yields the general solution
\[u(x,y)=F(y)(-1+x)^{-y},\]
which together with $u(0,y)=1$ gives the answer $u(x,y)=(1-x)^{-y}$, as is well known.
\item \textit{Stirling numbers of the second kind}. Here $(\alpha,\beta,\gamma,\alpha',\beta',\gamp )=(0,1,0,0,0,1)$, and the Maple package finds the general solution $u(x,y)=F(ye^x)e^{-y}$. The condition $u(0,y)=1$ chooses the answer \[u(x,y)=\exp{(y(e^x-1))}.\]
 \end{itemize}
\subsection{A more general case.} Take the case in which $\alpha'=0$ but otherwise the remaining five parameters are unrestricted. Now the \textit{pdsolve} command gets the general solution, and after applying the initial data $u(0,y)=1$ the generating function for this case is obtained in the form
\begin{equation}
\label{eq:uform}
u(x,y)=\frac{(1-\alpha x)^{-\frac{\gamma}{\alpha}}}{\left(1+\frac{\beta' }{\beta}y(1-(1-\alpha x)^{-\frac{\beta}{\alpha }})\right)^{1+\frac{\gamp }{\beta' }}}.
\end{equation}
By expanding this in a power series we find a very explicit form for the generalized binomial coefficients, viz.
\begin{equation}
\label{eq:formula}
{n\bars k}=\frac{\alpha^n}{k!}\left(\frac{\beta'}{\beta}\right)^{\kern-3pt k}\left(\frac{\gamp }{\beta'}+1\right)^{\kern-3pt\overline{k}}\ \sum_{j=0}^k(-1)^{k-j}{k\choose j}\left(\frac{\beta j+\gamma}{\alpha}\right)^{\kern-2pt \overline{n}}
\end{equation}
in which $x^{\overline{r}}$ is the usual rising factorial. A close inspection of this result shows that it is the product of a polynomial in $\alpha, \beta, \gamma$ only times a polynomial in $\beta',\gamp $ only. We can make this more explicit by expanding  the rising factorials into powers. This is a tedious calculation, but the result is that
\begin{eqnarray*}
{n\bars k}&=&\sum_{\ell,r,m}{k\brack \ell}{n\brack r}{r\choose m}{m\brace k}\alpha^{n-r}\beta^{m-k}\gamma^{r-m}(\beta')^{k-\ell}(\gamp )^{\ell}\\
&=&\sum_{\scriptstyle {{i_1+i_2+i_3=n-k}\atop {j_1+j_2=k}}}{k\brack j_2}{n\brack n-i_1}{n-i_1\choose k+i_2}{k+i_2\brace k}\alpha^{i_1}\beta^{i_2}\gamma^{i_3}(\beta')^{j_1}(\gamp )^{j_2}
\end{eqnarray*}
The factorization of the solution into a function of $\alpha,\beta,\gamma$ only, times a function of $\beta'$, $\gamp $ only, is
\[{n\bars k}=\left(\sum_{i_1+i_2+i_3=n-k}{n\brack n-i_1}{n-i_1\choose k+i_2}{k+i_2\brace k}\alpha^{i_1}\beta^{i_2}\gamma^{i_3}\right)\left(\prod_{j=1}^k(\gamp +j\beta')\right).\]
The appearance of three kinds of delimited Pascal-triangle sequences, viz. both kinds of Stirling numbers and the binomial coefficients, in the first factor is striking.

As an example we have, when $\alpha'=0$ and the other five parameters are unrestricted,
\[{4\bars 2}=(11\alpha^2+18\alpha\beta+7\beta^2+18\alpha\gamma+12\beta\gamma+6\gamma^2)(\beta'+\gamp )(2\beta'+\gamp ).\]
Additionally, we discover a major difference between this case where $\alpha'=0$ and the general case in that the general case does not exhibit such a factorization as above. For example, if all six of the parameters are unrestricted, we find that
\begin{eqnarray*}{4\bars 2}&=&
18\,{\alpha}\,{\gamma}\,{{\alpha'}}^2 + 7\,{\beta}\,{\gamma}\,{{\alpha'}}^2 + 11\,{{\gamma}}^2\,{{\alpha'}}^2 +
  18\,{{\alpha}}^2\,{\alpha'}\,{\beta'} + 29\,{\alpha}\,{\beta}\,{\alpha'}\,{\beta'} +
  11\,{{\beta}}^2\,{\alpha'}\,{\beta'}\\
&& + 51\,{\alpha}\,{\gamma}\,{\alpha'}\,{\beta'} +
  29\,{\beta}\,{\gamma}\,{\alpha'}\,{\beta'} + 22\,{{\gamma}}^2\,{\alpha'}\,{\beta'} +
  22\,{{\alpha}}^2\,{{\beta'}}^2 + 36\,{\alpha}\,{\beta}\,{{\beta'}}^2 + 14\,{{\beta}}^2\,{{\beta'}}^2\\
&&+
  36\,{\alpha}\,{\gamma}\,{{\beta'}}^2 + 24\,{\beta}\,{\gamma}\,{{\beta'}}^2 + 12\,{{\gamma}}^2\,{{\beta'}}^2 +
  18\,{{\alpha}}^2\,{\alpha'}\,{\gamp } + 29\,{\alpha}\,{\beta}\,{\alpha'}\,{\gamp } +
  11\,{{\beta}}^2\,{\alpha'}\,{\gamp }\\
&& + 44\,{\alpha}\,{\gamma}\,{\alpha'}\,{\gamp } +
  26\,{\beta}\,{\gamma}\,{\alpha'}\,{\gamp } + 18\,{{\gamma}}^2\,{\alpha'}\,{\gamp } +
  33\,{{\alpha}}^2\,{\beta'}\,{\gamp } + 54\,{\alpha}\,{\beta}\,{\beta'}\,{\gamp } +
  21\,{{\beta}}^2\,{\beta'}\,{\gamp }\\
&&+ 54\,{\alpha}\,{\gamma}\,{\beta'}\,{\gamp } +
  36\,{\beta}\,{\gamma}\,{\beta'}\,{\gamp } + 18\,{{\gamma}}^2\,{\beta'}\,{\gamp } +
  11\,{{\alpha}}^2\,{{\gamp }}^2 + 18\,{\alpha}\,{\beta}\,{{\gamp }}^2 + 7\,{{\beta}}^2\,{{\gamp }}^2\\
&& +
  18\,{\alpha}\,{\gamma}\,{{\gamp }}^2 + 12\,{\beta}\,{\gamma}\,{{\gamp }}^2 + 6\,{{\gamma}}^2\,{{\gamp }}^2,
\end{eqnarray*}
which does not factor as above.

Next let's get the solution to (\ref{eq:pde}), together with the initial condition $u(0,y)=1$, for a completely general set of parameters $(\alpha,\beta,\gamma,\alpha',\beta',\gamp )$. The package \texttt{pdsolve} first returns the general solution $u(x,y)$ in the form
\begin{equation}\label{eq:gencase}
(\beta +\beta \,'y)^{\frac{\gamma }{\beta }-1-\frac{\gamma '}{\beta \,'}}y^{-\frac{\gamma }{\beta }}f\left(\frac{1}{\alpha }y^{-\frac{\alpha }{\beta }}\left((\beta +\beta \,'y)^{\frac{\alpha }{\beta }-\frac{\alpha '}{\beta \,'}}\alpha x-\hypf{2}{1}{-\frac{\alpha }{\beta },1+\frac{\alpha '}{\beta \,'}-\frac{\alpha }{\beta }}{1-\frac{\alpha }{\beta }}{-y\frac{\beta \,'}{\beta }}\beta ^{\frac{\alpha }{\beta }-\frac{\alpha '}{\beta \,'}}\right)\right)
\end{equation}
The condition $u(0,y)=1$ gives
\[f\left(-\frac{1}{\alpha }y^{-\frac{\alpha }{\beta }}\hypf{2}{1}{-\frac{\alpha }{\beta },1+\frac{\alpha '}{\beta \,'}-\frac{\alpha }{\beta }}{1-\frac{\alpha }{\beta }}{-y\frac{\beta \,'}{\beta }}\beta ^{\frac{\alpha }{\beta }-\frac{\alpha '}{\beta \,'}}\right)=(\beta +\beta \,'y)^{-\frac{\gamma }{\beta }+1+\frac{\gamma '}{\beta \,'}}y^{\frac{\gamma }{\beta }}.\]
Hence let $y=y(t)$ be the inverse function of the argument of $f$ in the above. If we replace $y$ by $y(t)$ throughout, we obtain
\[f(t)=(\beta +\beta \,'y(t))^{-\frac{\gamma }{\beta }+1+\frac{\gamma '}{\beta \,'}}y(t)^{\frac{\gamma }{\beta }}.\]
With this choice of $f$, the general solution is given by (\ref{eq:gencase}), though in quite an unwieldy form.

\subsection{Reality of the zeros}
We will now study the polynomials
\[\phi_n(x)=\sum_k{n\bars k}x^k\qquad (n=0,1,2,\dots),\]
with a view to describing when they will have real zeros only. In the cases of binomial coefficients, both kinds of Stirling numbers, and others, this is true. We limit the discussion to the case where all of $\alpha,\beta,\gamma,\alpha',\beta',\gamp $ are nonnegative. Even in this case, however, it is not true that the polynomials must have real zeros only. If we take
\[(\alpha,\beta,\gamma,\alpha',\beta',\gamp )=(3,1,1,0,1,0),\]
for example, then $\phi_3(x)=2(x+1)(3x^2+12x+14)$, which has two nonreal zeros.

If we multiply the ``Pascal triangle'' recurrence (\ref{eq:recur}) by $x^k$ and sum on $k$ we find the recurrence for the $\phi_n$ polynomials as
\begin{equation}
\label{eq:phire}
\phi_{n+1}(x)=((\alpha n+\gamma)+(\alpha'n+\beta'+\gamp )x)\phi_n(x)+(\beta+\beta' x)x\phi_n'(x).
\end{equation}

First suppose that $\beta=\beta'=0$. Then (\ref{eq:phire}) gives the explicit expression
\[\phi_n(x)=\prod_{j=0}^{n-1}((\alpha j+\gamma)+(\alpha' j+\gamp )x),\]
which plainly has real zeros only.

Next, suppose $\beta'=0$, and $\beta, \gamma, \gamp >0$.  Let, for $n$ fixed,
\[u=\alpha (n-1)+\gamma,\quad v=\alpha'(n-1)+\gamp .\]
In this case the $\phi_n$ recurrence can be cast in the form
\begin{equation}
\label{eq:indra}
\phi_{n}(x)=
\begin{cases}
\beta x\, |x|^{-u/\beta}e^{-vx/\beta}\frac{d}{dx}\left\{|x|^{u/\beta}e^{vx/\beta}\phi_{n-1}(x)\right\},&\mathrm{if}\ \text{$x\neq 0$, $n\ge 1$};\\
\prod_{j=0}^{n-1}(\gamma+j\alpha),&\mathrm{if}\ \text{$x=0$},
\end{cases}
\end{equation}
with $\phi_0(x)=1$. Now $u,v,\beta$ are positive. If, inductively, $\phi_{n-1}(x)$ has $n-1$ negative real zeros then the function $|x|^{u/\beta}e^{vx/\beta}\phi_{n-1}(x)$ vanishes at those $n-1$ points and at $-\infty$ and at the origin. By Rolle's theorem, the derivative of this function must vanish at $n$ points, whence $\phi_n(x)$ has $n$ negative real zeros.

This shows that if $\beta'=0$ and $\alpha,\alpha',\beta$ are nonnegative, and  $\gamma$ and $\gamp $ are positive, then the zeros are all real and negative. By an obvious continuity argument, it must be that if $\beta'=0$ and the remaining parameters are all nonnegative then the zeros are real and nonpositive. Hence we have,
\begin{theorem}
If $\beta'=0$ and $\alpha, \beta, \gamma, \alpha',\gamp $  are nonnegative, then all of the polynomials $\{\phi_n(x)\}_0^{\infty}$ have real zeros only. There are examples with $\beta'\neq 0$ in which non-real zeros exist.
\end{theorem}
\subsection{A question: the row sums}
We define the row sums of the generalized Pascal triangle,
\begin{equation}
\label{eq:rowsm}
\rho(n)=\sum_k{n\bars k}\qquad (n=0,1,2,\dots).
\end{equation}
In the classical cases these row sums are variously $2^n$, $n!$, Bell$(n)$, etc. We ask for information about their asymptotic growth rate,  under the assumption that $\alpha'=0$ so that the explicit generating function (\ref{eq:uform}) is available. Indeed from (\ref{eq:uform}) we have, by putting $y:=1$, the generating function of the row sums as
\[F(x)\eqdef \sum_n\rho(n)\frac{x^n}{n!}=\frac{(1-\alpha x)^{-\frac{\gamma}{\alpha}}}{\left(1+\frac{\beta' }{\beta}(1-(1-\alpha x)^{-\frac{\beta}{\alpha }})\right)^{1+\frac{\gamp }{\beta' }}}.\]
From this, the exact formula for the row sums is
\[\rho(n)=\alpha^nn!\sum_{j\ge 0}{j+\frac{\gamp }{\beta'}\choose j}{n+\frac{j\beta+\gamma}{\alpha}-1\choose n}\left(1+\frac{\beta'}{\beta}\right)^{-j-1-\gamp /\beta'}.\]
What is the asymptotic behavior of these row sums for large $n$? This will depend on the analytic character of $F(x)$ above, and will no doubt require a careful case analysis.

\end{document}